\documentclass[review]{elsarticle}
\usepackage{lineno,hyperref}
\usepackage{amssymb}
\usepackage{amsfonts}
\usepackage{amsmath}
\usepackage{amsthm}
 \newtheorem{theorem}{Theorem}[section]
 \newtheorem{corollary}[theorem]{Corollary}

 \newtheorem{definition}[theorem]{Definition}
 \newtheorem{remark}[theorem]{Remark}
 \newtheorem{example}[theorem]{Example}
\modulolinenumbers[5]
\newenvironment{bsmallmatrix}
  {\left(\begin{smallmatrix}}
  {\end{smallmatrix}\right)}

\begin{document}

\begin{frontmatter}

\title{Multivariable Conformable Fractional Calculus}


\author[mymainaddress]{Nazl{\i} Yaz{\i}c{\i} G\"oz\"utok }
\ead{nazliyazici@ktu.edu.tr}

\author[mymainaddress]{U\u{g}ur  G\"oz\"utok\corref{mycorrespondingauthor}}
\cortext[mycorrespondingauthor]{Corresponding author}
\ead{ugurgozutok@ktu.edu.tr}

\address[mymainaddress]{Department of Mathematics \\ 
Karadeniz Technical University \\
61080, Trabzon, TURKEY}

\begin{abstract}
Conformable fractional derivative is introduced by the authors Khalil et al. In this study we develop their concept and introduce multivariable conformable derivative for a vector valued function with several variables.
\end{abstract}

\begin{keyword}
Conformable fractional derivative\sep multivariable conformable calculus\sep conformable partial derivatives\sep conformable jacobian matrix 
\MSC[2010] 26B12\sep 26A33
\end{keyword}

\end{frontmatter}

\section{Introduction}

For many years, many definitions of fractional derivative have been introduced by various researchers. One of them is the Riemann-Liouville fractional derivative and the second one is the so-called Caputo derivative. But these are not all purpose definitions. Recently, a new fractional derivative has been introduced in \citep{1} and one can see that the new derivative suggested in these papers satisfies all the properties of the standart one. However, definitions given in the literature are only for the real valued functions. In this paper, we give conformable fractional derivative definition for the vector valued functions with several real variables. Our paper has the following order: In Section 2, some basic definitions and theorems appeared in the literature are given. In Section 3, $\alpha -$derivative of a vector valued function, conformable Jacobian matrix are defined; relation between $\alpha -$derivative and usual derivative of a vector valued function is revealed; chain rule for multivariable conformable derivative is given. In Section 4, conformable partial derivatives of a real valued function with $n-$variables is defined and relation between conformable Jacobian matrix and conformable partial derivatives is given.       

\section{Basic definitions and theorems}
In this section we will give some definitions and properties introduced in \citep{1,2}

\begin{definition}\label{Def2.1}
Given a function $f:[0,\infty)\longrightarrow \mathbb{R}$. The conformable derivative of the function $f$ of order $\alpha$ is defined by
\begin{equation}\label{eq1}
T_{\alpha}(f)(x)=\lim_{h\to 0} \dfrac{f(x+hx^{1-\alpha})-f(x)}{h}
\end{equation}
for all $x>0$, $\alpha\in (0,1)$.
\end{definition}
\begin{theorem}\label{Theo2.2}
If a function $f:[0,\infty)\longrightarrow \mathbb{R}$ is $\alpha -$differentiable at $t_0 >0$, $\alpha\in (0,1]$, then $f$ is continuous at $t_0$.
\end{theorem}
\begin{theorem}\label{Theo2.3}
Let $\alpha\in (0,1]$ and $f,g$ be $\alpha -$differentiable at a point $t>0$. Then
\begin{itemize}
\item[(1)] $T_\alpha (af+bg)=aT_\alpha (f)+bT_\alpha (g)$, for all $a,b\in\mathbb{R}$.
\item[(2)] $T_\alpha (t^p)=pt^{p-\alpha}$ for all $p\in\mathbb{R}$.
\item[(3)] $T_\alpha (\lambda) =0$, for all constant functions $f(t)=\lambda$.
\item[(4)] $T_\alpha (fg)=fT_\alpha (g)+gT_\alpha (f)$.
\item[(5)] $T_\alpha (\dfrac{f}{g})=\dfrac{gT_\alpha (f)-fT_\alpha (g)}{g^2}$.
\item[(6)] If, in addition, f is differentiable, then $T_\alpha (f)(t)=t^{1-\alpha}\dfrac{d}{dt}f(t)$.
\end{itemize}
\end{theorem}
\begin{theorem}\label{Theo2.4}
Assume $f,g:(0,\infty)\longrightarrow\mathbb{R}$ be two $\alpha -$differentiable functions where $\alpha\in (0,1]$. Then $g\circ f$ is $\alpha -$differentiable and for all $t$ with $t\neq 0$  and $f(t)\neq 0$ we have
\begin{equation}\label{eq2}
T_\alpha (g\circ f)(t)=T_\alpha (g)(f(t))T_\alpha (f)(t)f(t)^{\alpha -1}.
\end{equation} 
\end{theorem}

\section{\texorpdfstring{$\alpha -$}{Lg}Derivative of a Vector Valued Function}
\begin{definition}\label{Def3.1}
Let $f$ be a vector valued function with $n$ real variables such that $f(x_1,...,x_n)=(f_1(x_1,...,x_n),...,f_m(x_1,...,x_n))$. Then we say that $f$ is $\alpha -$differentiable at $a=(a_1,...,a_n)\in\mathbb{R}^n$ where each $a_i>0$, if there is a linear transformation $L:\mathbb{R}^n\longrightarrow\mathbb{R}^m$ such that
\begin{equation}\label{eq3}
\lim_{h\to 0} \dfrac{\lVert f(a_1+h_1a_1^{1-\alpha},...,a_n+h_na_n^{1-\alpha})-f(a_1,...,a_n)-L(h) \rVert}{\lVert h \rVert}=0
\end{equation}
where $h=(h_1,...,h_n)$ and $\alpha\in (0,1]$. The linear transformation is denoted by $D^\alpha f(a)$ and called the conformable derivative of $f$ of order $\alpha$ at $a$.    
\end{definition}
\begin{remark}\label{rem3.2}
For $m=n=1$, Definition \ref{Def3.1} equivalent to Definition \ref{Def2.1}. 
\end{remark}
\begin{theorem}\label{theo3.3}
Let $f$ be a vector valued function with $n$ variables. If $f$ is $\alpha -$differentiable at $a=(a_1,...,a_n)\in\mathbb{R}^n$, each $a_i>0$, then there is a unique linear transformation $L:\mathbb{R}^n\longrightarrow\mathbb{R}^m$ such that
\begin{equation*}
\lim_{h\to 0} \dfrac{\lVert f(a_1+h_1a_1^{1-\alpha},...,a_n+h_na_n^{1-\alpha})-f(a_1,...,a_n)-L(h) \rVert}{\lVert h \rVert}=0.
\end{equation*}
\end{theorem}
\begin{proof}
Suppose $M:\mathbb{R}^n\longrightarrow\mathbb{R}^m$ satisfies
\begin{equation*}
\lim_{h\to 0} \dfrac{\lVert f(a_1+h_1a_1^{1-\alpha},...,a_n+h_na_n^{1-\alpha})-f(a_1,...,a_n)-M(h) \rVert}{\lVert h \rVert}=0.
\end{equation*}
Then,
\begin{multline*}
\lim_{h\to 0}\dfrac{\lVert L(h)-M(h) \rVert}{\lVert h \rVert} \\
\leq\lim_{h\to 0} \dfrac{\lVert L(h)-(f(a_1+h_1a_1^{1-\alpha},...,a_n+h_na_n^{1-\alpha})-f(a)) \rVert}{\lVert h \rVert} \\
+\lim_{h\to 0} \dfrac{\lVert (f(a_1+h_1a_1^{1-\alpha},...,a_n+h_na_n^{1-\alpha})-f(a))-M(h) \rVert}{\lVert h \rVert}=0
\end{multline*}
If $x\in\mathbb{R}^n$, then $tx\rightarrow 0$ as $t\rightarrow 0$. Hence for $x\neq 0$ we have
\begin{align*}
0=\lim_{h\to 0}\dfrac{\lVert L(tx)-M(tx) \rVert}{\lVert tx \rVert}=\dfrac{\lVert L(x)-M(x) \rVert}{\lVert x \rVert}.
\end{align*}
Therefore $L(x)=M(x).$
\end{proof}
\begin{example}\label{exm3.4}
Let us consider the function $f$ defined by $f(x,y)=sinx$ and the point $(a,b)\in\mathbb{R}^2$ such that $a,b>0$, then $D^\alpha f(a,b)=L$ satisfies $L(x,y)=xa^{1-\alpha}\cos a$.

To prove this, note that
\begin{align*}
& \lim_{(h_1,h_2)\to (0,0)}\dfrac{\mid f(a+h_1a^{1-\alpha},b+h_2b^{1-\alpha})-f(a,b)-L(h_1,h_2)\mid}{\lVert (h_1,h_2)\rVert} \\
& =\lim_{(h_1,h_2)\to (0,0)}\dfrac{\mid \sin (a+h_1a^{1-\alpha})-\sin a-h_1a^{1-\alpha}\cos a\mid}{\sqrt{h_1^2+h_2^2}} \\
& \leq \lim_{h_1\to 0}\dfrac{\mid \sin (a+h_1a^{1-\alpha})-\sin a-h_1a^{1-\alpha}\cos a\mid}{\mid h_1 \mid} \\
& =\lim_{h_1\to 0}\mid \dfrac{\sin (a+h_1a^{1-\alpha})-\sin a}{h_1}-a^{1-\alpha}\cos a \mid \\
& =\mid a^{1-\alpha}\cos a -a^{1-\alpha}\cos a \mid =0
\end{align*}
\end{example}
\begin{definition}\label{Def3.5}
Consider the matrix of the linear transformation $D^\alpha f(a):\mathbb{R}^n\longrightarrow\mathbb{R}^m$ with respect to the usual bases of $\mathbb{R}^n$ and $\mathbb{R}^m$. This $m\times n$ matrix is called the conformable Jacobian matrix of $f$ at $a$, and denoted by $f^\alpha (a)$.
\end{definition}
\begin{example}\label{exm3.6}
If $f(x,y)=\sin x$, then $f^\alpha (a,b)=\begin{bmatrix}
a^{1-\alpha}\cos a & 0
\end{bmatrix}$
\end{example}
\begin{theorem}\label{Theo3.7}
Let $f$ be $\alpha -$differentiable at $a=(a_1,...,a_n)\in\mathbb{R}^n$, each $a_i>0$. If $f$ is differentiable at $a$, then
\begin{align*}
D^\alpha f(a)=Df(a)\circ L_a^{1-\alpha}
\end{align*}
where $Df(a)$ is the usual derivative of $f$ and $L_a^{1-\alpha}$ is the linear transformation from $\mathbb{R}^n$ to $\mathbb{R}^n$ defined by $L_a^{1-\alpha}(x_1,...,x_n)=(a_1^{1-\alpha}x_1,...,a_n^{1-\alpha}x_n)$.
\end{theorem}
\begin{proof}
It is suffices to show that
\begin{align*}
\lim_{h\to 0} \dfrac{\lVert f(a_1+h_1a_1^{1-\alpha},...,a_n+h_na_n^{1-\alpha})-f(a_1,...,a_n)-Df(a)\circ L_a^{1-\alpha}(h) \rVert}{\lVert h \rVert}=0.
\end{align*}
Let $\epsilon =(\epsilon _1,...,\epsilon _n)=(h_1a_1^{1-\alpha},...,h_na_n^{1-\alpha})$, then $\epsilon \rightarrow 0$ as $h\rightarrow 0$. On the other hand, if we put $M=\max\{ (a_i^{1-\alpha})^2\mid a_i>0, i=1,2,...,n\}>0$, then 
\begin{align*}
\lVert \epsilon \rVert =\sqrt{h_1^2(a_1^{1-\alpha})^2+...+h_n^2(a_n^{1-\alpha})^2}\leq \sqrt{h_1^2M+...+h_n^2M}=\sqrt{nM}\lVert h \rVert .
\end{align*}
Hence we have
\begin{align*}
\dfrac{1}{\sqrt{nM}} \lVert \epsilon \rVert \leq \lVert h \rVert .
\end{align*}
Finally,
\begin{align*}
& \lim_{h\to 0} \dfrac{\lVert f(a_1+h_1a_1^{1-\alpha},...,a_n+h_na_n^{1-\alpha})-f(a_1,...,a_n)-Df(a)\circ L_a^{1-\alpha}(h) \rVert}{\lVert h \rVert} \\
& =\lim_{h\to 0} \dfrac{\lVert f(a_1+h_1a_1^{1-\alpha},...,a_n+h_na_n^{1-\alpha})-f(a)-Df(a)(h_1a_1^{1-\alpha},...,h_na_n^{1-\alpha}) \rVert}{\lVert h \rVert} \\
& \leq \lim_{\epsilon\to 0}\frac{\lVert f(a+\epsilon)-f(a)-Df(a)(\epsilon) \rVert}{\dfrac{1}{\sqrt{nM}} \lVert \epsilon \rVert} \\
& =\sqrt{nM}\lim_{\epsilon\to 0}\frac{\lVert f(a+\epsilon)-f(a)-Df(a)(\epsilon) \rVert}{ \lVert \epsilon \rVert}=\sqrt{nM}.0=0.
\end{align*}
This completes the proof. 
\end{proof}
Theorem \ref{Theo3.7} is the generalized case of the part (6) of Theorem \ref{Theo2.3}. Also matrix form of the Theorem \ref{Theo3.7} is given by the following:
\begin{equation*}
f^\alpha (a)=f'(a).\begin{bsmallmatrix}
a_1^{1-\alpha} &      0         & \hdots  & 0      \\
0              & a_1^{1-\alpha} & \hdots  & 0      \\
\vdots         & \vdots         & \ddots  & \vdots \\      
0              &      0         & \hdots        & a_n^{1-\alpha} 
\end{bsmallmatrix},
\end{equation*}
where $f'(a)$ is the usual Jacobian of $f$ and $\begin{bsmallmatrix}
a_1^{1-\alpha} &      0         & \hdots  & 0      \\
0              & a_1^{1-\alpha} & \hdots  & 0      \\
\vdots         & \vdots         & \ddots  & \vdots \\      
0              &      0         & \hdots        & a_n^{1-\alpha} 
\end{bsmallmatrix}$ is the matrix corresponding to linear transformation $L_a^{1-\alpha}$.
\begin{theorem}\label{Theo3.8}
If a vector valued function $f$ with $n$ variables is $\alpha -$differentiable at $a=(a_1,...,a_n)\in \mathbb{R}^n$, each $a_i>0$, then $f$ is continuous at $a\in \mathbb{R}^n$.
\end{theorem}
\begin{proof}
Since 
\begin{align*}
&{\lVert f(a_1+h_1a_1^{1-\alpha},...,a_n+h_na_n^{1-\alpha})-f(a_1,...,a_n)\rVert}\\
&=\dfrac{\lVert f(a_1+h_1a_1^{1-\alpha},...,a_n+h_na_n^{1-\alpha})-f(a)-L(h)+L(h)\rVert}{\lVert h \rVert}{\lVert h \rVert} \\
&\leq \dfrac{\lVert f(a_1+h_1a_1^{1-\alpha},...,a_n+h_na_n^{1-\alpha})-f(a)-L(h)\rVert}{\lVert h \rVert}{\lVert h \rVert}+{\lVert L(h) \rVert}.
\end{align*}
We have
\begin{align*}
&{\lVert f(a_1+h_1a_1^{1-\alpha},...,a_n+h_na_n^{1-\alpha})-f(a_1,...,a_n)\rVert}\\
&\leq \dfrac{\lVert f(a_1+h_1a_1^{1-\alpha},...,a_n+h_na_n^{1-\alpha})-f(a)-L(h)\rVert}{\lVert h \rVert}{\lVert h \rVert}+{\lVert L(h) \rVert}.
\end{align*}
By taking limits of the two sides of the inequality as $ h\rightarrow 0$, we have 
\begin{align*}
& \lim_{h\to 0} {\lVert f(a_1+h_1a_1^{1-\alpha},...,a_n+h_na_n^{1-\alpha})-f(a_1,...,a_n) \rVert} \\
&\leq \lim_{h\to 0} \dfrac{\lVert f(a_1+h_1a_1^{1-\alpha},...,a_n+h_na_n^{1-\alpha})-f(a)-L(h)\rVert}{\lVert h \rVert}\lim_{h\to 0} {\lVert h \rVert} +\lim_{h\to 0} {\lVert L(h) \rVert}\\
&=0.
\end{align*}
Let $(\epsilon_1,...,\epsilon_n)=(h_1a_1^{1-\alpha},...,h_na_n^{1-\alpha})$, then $ \epsilon\rightarrow 0$ as $ h\rightarrow 0$.
Since
$$\lim_{\epsilon\to 0} {\lVert f(a+\epsilon)-f(a) \rVert} \leq 0$$
we have $$\lim_{\epsilon\to 0} {\lVert f(a+\epsilon)-f(a) \rVert} = 0.$$
Hence $f$ is continuous at $a\in \mathbb{R}^n$.
\end{proof}
\begin{theorem}\label{Theo3.9}
(Chain Rule) Let $x\in\mathbb{R}^n$, $y\in\mathbb{R}^m$. If $f(x)=(f_1(x),...,f_m(x))$ is $\alpha -$differentiable at $a=(a_1,...,a_n)\in\mathbb{R}^n$, each $a_i>0$ such that $\alpha\in (0,1]$ and $g(y)=(g_1(y),...,g_p(y))$ is $\alpha -$differentiable at $f(a)\in\mathbb{R}^m$, all $f_i(a)>0$ such that $\alpha\in (0,1]$. Then the composition $g\circ f$ is $\alpha -$differentiable at $a$ and
\begin{equation}
D^\alpha (g\circ f)(a)=D^\alpha g(f(a))\circ f(a)^{\alpha -1} \circ D^\alpha f(a)
\end{equation}
where $f(a)^{\alpha -1}$ is the linear transformation from $\mathbb{R}^m$ to $\mathbb{R}^m$ defined by \\ $f(a)^{\alpha -1}(x_1,...,x_m)=(x_1f_1(a)^{\alpha -1},...,x_mf_m(a)^{\alpha -1}).$
\end{theorem}
\begin{proof}
Let $L=D^\alpha f(a)$, $M=D^\alpha g(f(a))$. If we define
\begin{itemize}
\item[(i)]$\varphi (a_1+h_1a_1^{1-\alpha},...,a_n+h_na_n^{1-\alpha}) \\
=f(a_1+h_1a_1^{1-\alpha},...,a_n+h_na_n^{1-\alpha})-f(a)-L(h)$,
\item[(ii)]$\psi (f_1(a)+k_1f_1(a)^{1-\alpha},...,f_n(a)+k_nf_n(a)^{1-\alpha}) \\
=g(f_1(a)+k_1f_1(a)^{1-\alpha},...,f_n(a)+k_nf_n(a)^{1-\alpha})-g(f(a))-M(k)$,
\item[(iii)]$\rho (a_1+h_1a_1^{1-\alpha},...,a_n+h_na_n^{1-\alpha}) \\
=g\circ f(a_1+h_1a_1^{1-\alpha},...,a_n+h_na_n^{1-\alpha})-g\circ f(a)-M\circ f(a)^{\alpha -1}\circ L(h)$,
\end{itemize}
then
\begin{itemize}
\item[(iv)]$\lim_{h\to 0}\dfrac{\lVert \varphi (a_1+h_1a_1^{1-\alpha},...,a_n+h_na_n^{1-\alpha}) \rVert}{\lVert h \rVert}=0$,
\item[(v)]$\lim_{k\to 0}\dfrac{\lVert \psi (f_1(a)+k_1f_1(a)^{1-\alpha},...,f_n(a)+k_nf_n(a)^{1-\alpha}) \rVert}{\lVert k \rVert}=0$,
\end{itemize}
and we must show that
\begin{equation*}
lim_{h\to 0}\dfrac{\lVert \rho (a_1+h_1a_1^{1-\alpha},...,a_n+h_na_n^{1-\alpha})\rVert}{\lVert h \rVert}=0
\end{equation*}
Now,
\begin{align*}
& \rho (a_1+h_1a_1^{1-\alpha},...,a_n+h_na_n^{1-\alpha}) \\
& =g(f(a_1+h_1a_1^{1-\alpha},...,a_n+h_na_n^{1-\alpha}))-g(f(a))-M\circ f(a)^{\alpha -1}\circ L(h) \\
& =g(f_1(a_1+h_1a_1^{1-\alpha},...,a_n+h_na_n^{1-\alpha}),..., f_m(a_1+h_1a_1^{1-\alpha},...,a_n+h_na_n^{1-\alpha})) \\
& -g(f(a))-M\circ f(a)^{\alpha -1}(f(a_1+h_1a_1^{1-\alpha},...,a_n+h_na_n^{1-\alpha})-f(a) \\
& -\varphi (a_1+h_1a_1^{1-\alpha},...,a_n+h_na_n^{1-\alpha})) \qquad by (i) \\
& =[g(f_1(a_1+h_1a_1^{1-\alpha},...,a_n+h_na_n^{1-\alpha}),..., f_m(a_1+h_1a_1^{1-\alpha},...,a_n+h_na_n^{1-\alpha})) \\
& -g(f(a))-M(f(a)^{\alpha -1}(f(a_1+h_1a_1^{1-\alpha},...,a_n+h_na_n^{1-\alpha})-f(a)))] \\
& +M\circ f(a)^{\alpha -1}(\varphi (a_1+h_1a_1^{1-\alpha},...,a_n+h_na_n^{1-\alpha})) \\
& =[g(f_1(a_1+h_1a_1^{1-\alpha},...,a_n+h_na_n^{1-\alpha}),..., f_m(a_1+h_1a_1^{1-\alpha},...,a_n+h_na_n^{1-\alpha})) \\
& -g(f(a))-M(f(a)^{\alpha -1}(f_1(a_1+h_1a_1^{1-\alpha},...,a_n+h_na_n^{1-\alpha})-f_1(a),..., \\
& f_m(a_1+h_1a_1^{1-\alpha},...,a_n+h_na_n^{1-\alpha})-f_m(a)))] \\
& +M\circ f(a)^{\alpha -1}(\varphi (a_1+h_1a_1^{1-\alpha},...,a_n+h_na_n^{1-\alpha})) 
\end{align*}
\begin{align*}
& =[g(f_1(a_1+h_1a_1^{1-\alpha},...,a_n+h_na_n^{1-\alpha}),..., f_m(a_1+h_1a_1^{1-\alpha},...,a_n+h_na_n^{1-\alpha})) \\
& -g(f(a))-M([f_1(a_1+h_1a_1^{1-\alpha},...,a_n+h_na_n^{1-\alpha})-f_1(a)]f_1(a)^{\alpha -1},..., \\
& [f_m(a_1+h_1a_1^{1-\alpha},...,a_n+h_na_n^{1-\alpha})-f_m(a)]f_m(a)^{\alpha -1})] \\
& +M\circ f(a)^{\alpha -1}(\varphi (a_1+h_1a_1^{1-\alpha},...,a_n+h_na_n^{1-\alpha}))
\end{align*}
If we put $u_i=[f_i(a_1+h_1a_1^{1-\alpha},...,a_n+h_na_n^{1-\alpha})-f_i(a)]f_i(a)^{\alpha -1}$, $i=1,...,n$, then we have $f_i(a_1+h_1a_1^{1-\alpha},...,a_n+h_na_n^{1-\alpha})=f_i(a)+u_if_i(a)^{1-\alpha}$ and $u\rightarrow 0$ as $h\rightarrow 0$.
Therefore,
\begin{align*}
& \rho (a_1+h_1a_1^{1-\alpha},...,a_n+h_na_n^{1-\alpha}) \\
& =[g(f_1(a_1+h_1a_1^{1-\alpha},...,a_n+h_na_n^{1-\alpha}),..., f_m(a_1+h_1a_1^{1-\alpha},...,a_n+h_na_n^{1-\alpha})) \\
& -g(f(a))-M(u)]+M\circ f(a)^{\alpha -1}(\varphi (a_1+h_1a_1^{1-\alpha},...,a_n+h_na_n^{1-\alpha})) \\
& =\psi(f_1(a)+u_1f_1(a)^{1-\alpha},...,f_m(a)+u_mf_m(a)^{1-\alpha}) \\
& +M\circ f(a)^{\alpha -1}(\varphi (a_1+h_1a_1^{1-\alpha},...,a_n+h_na_n^{1-\alpha})) \qquad by (ii).
\end{align*}
Thus we must show
\begin{itemize}
\item[(vi)]$\lim_{u\to 0}\dfrac{\lVert \psi(f_1(a)+u_1f_1(a)^{1-\alpha},...,f_m(a)+u_mf_m(a)^{1-\alpha}) \rVert}{\lVert u \rVert}=0,$
\item[(vii)]$\lim_{h\to 0}\dfrac{\lVert M\circ f(a)^{\alpha -1}(\varphi (a_1+h_1a_1^{1-\alpha},...,a_n+h_na_n^{1-\alpha})) \rVert}{\lVert h \rVert}=0.$
\end{itemize}
For (vi), it is obvious from (v). For (vii), the linear transformation satisfies
\begin{align*}
\lVert M\circ f(a)^{\alpha -1}(\varphi (a_1+h_1a_1^{1-\alpha},...,a_n+h_na_n^{1-\alpha})) \rVert \\
\leq K\lVert \varphi (a_1+h_1a_1^{1-\alpha},...,a_n+h_na_n^{1-\alpha})\rVert
\end{align*}
such that $K>0$. Hence, from (iv), (vii) holds.
\end{proof} 
\begin{corollary}\label{Cor3.10}
For, $m=n=p=1$, Theorem \ref{Theo3.9} states that
\begin{align*}
T_\alpha (g\circ f)(a)=T_\alpha g(f(a))T_\alpha f(a) f(a)^{\alpha -1}.
\end{align*}
\end{corollary}
Above corollary says that Theorem \ref{Theo3.9} is the generalized case of Theorem \ref{Theo2.4}.
\begin{corollary}\label{Cor3.11}
Let all conditions of Theorem \ref{Theo3.9} be satisfied. Then
\begin{align*}
(g\circ f)^{\alpha}(a)=g^{\alpha}(f(a))\begin{bsmallmatrix}
f_1(a)^{\alpha -1} &      0     & \hdots   &   0   \\
0                & f_2(a)^{\alpha -1} & \hdots  &  0 \\
\vdots           &  \vdots        & \ddots & \vdots  \\
0 & 0 & \hdots & f_m(a)^{\alpha -1}
\end{bsmallmatrix}
f^{\alpha}(a)
\end{align*}
where $\begin{bsmallmatrix}
f_1(a)^{\alpha -1} &      0     & \hdots   &   0   \\
0                & f_2(a)^{\alpha -1} & \hdots  &  0 \\
\vdots           &  \vdots        & \ddots & \vdots  \\
0 & 0 & \hdots & f_m(a)^{\alpha -1}
\end{bsmallmatrix}$ is the matrix corresponding to the linear transformation $f(a)^{\alpha -1}$.
\end{corollary}
\begin{theorem}\label{Theo3.12}
Let $f$ be a vector valued function with $n$ variables such that $f(x_1,...,x_n)=(f_1(x_1,...,x_n),...,f_m(x_1,...,x_n))$. Then $f$ is $\alpha -$differentiable at $a=(a_1,...,a_n)\in\mathbb{R}^n $, each $a_i>0$ if and only if each $f_i$ is, and
\begin{align*}
D^\alpha f(a)= (D^\alpha f_1(a),...,D^\alpha f_m(a)).
\end{align*}
\end{theorem}
\begin{proof}
If each $f_i$ is $\alpha -$differentiable at $a$ and $L=(D^\alpha f_1(a),...,D^\alpha f_m(a))$, then
\begin{align*}
&f(a_1+h_1a_1^{1-\alpha},...,a_n+h_na_n^{1-\alpha})-f(a)-L(h) \\
&=(f_1(a_1+h_1a_1^{1-\alpha},...,a_n+h_na_n^{1-\alpha})-f_1(a)-D^\alpha f_1(a)(h),..., \\
& f_m(a_1+h_1a_1^{1-\alpha},...,a_n+h_na_n^{1-\alpha})-f_m(a)-D^\alpha f_m(a)(h)).
\end{align*}
Therefore,
\begin{align*}
& \lim_{h\to 0}\dfrac{\lVert f(a_1+h_1a_1^{1-\alpha},...,a_n+h_na_n^{1-\alpha})-f(a)-L(h) \rVert}{\lVert h\rVert} \\
& \leq \lim_{h\to 0}\sum\limits_{i=1}^n\dfrac{\lVert f_i(a_1+h_1a_1^{1-\alpha},...,a_n+h_na_n^{1-\alpha})-f_i(a)-D^\alpha f_i(a)(h) \rVert}{\lVert h \rVert}=0.
\end{align*}
If, on the other hand, $f$ is $\alpha -$differentiable at $a$, then $f_i=\pi _if$ is $\alpha -$differentiable at $a$ by Theorem \ref{Theo3.9}.
\end{proof}
\begin{theorem}\label{Theo3.13}
Let $\alpha\in (0,1]$ and $f,g$ be $\alpha -$differentiable at a point $a=(a_1,...,a_n)\in\mathbb{R}^n$, each $a_i>0$. Then
\begin{itemize}
\item[(i)]$D^\alpha (\lambda f+\mu g)(a)=\lambda D^\alpha f(a)+\mu D^\alpha g(a)$ for all $\lambda ,\mu \in\mathbb{R}$.
\item[(ii)]$D^\alpha (fg)(a)=f(a)D^\alpha g(a)+g(a)D^\alpha f(a)$.
\end{itemize}
\end{theorem}
\begin{proof}
(i) follows from the definition, thus we omitted the proof of (i).

For (ii), let $a_1+h_1a_1^{1-\alpha},...,a_n+h_na_n^{1-\alpha}=A$, then
\begin{align*}
& \lim_{h\to 0}\dfrac{\lVert (fg)(A)-(fg)(a)-(f(a)D^\alpha g(a)+g(a)D^\alpha f(a))(h) \rVert}{\lVert h \rVert} \\
& \leq \lim_{h\to 0}\dfrac{\lVert f(A)g(A)-f(a)g(A)-g(A)D^\alpha f(a)(h)\rVert}{\lVert h \rVert} \\
& +\lim_{h\to 0}\dfrac{\lVert f(a)g(A)-f(a)g(a)-f(a)D^\alpha g(a)(h) \rVert}{\lVert h \rVert} \\
& +\lim_{h\to 0}\dfrac{\lVert g(A)D^\alpha f(a)(h)-g(a)D^\alpha f(a)(h)\rVert}{\lVert h\rVert} \\
& =\lim_{h\to 0}{\lVert g(A)\rVert} \dfrac{\lVert f(A)-f(a)-D^\alpha f(a)(h)\rVert}{\lVert h\rVert} \\
& +\lim_{h\to 0}{\lVert f(a)\rVert} \dfrac{\lVert g(A)-g(a)-D^\alpha g(a)(h)\rVert}{\lVert h\rVert}+\lim_{h\to 0}\lVert D^\alpha f(a)(h)\rVert\dfrac{\lVert g(A)-g(a)\rVert}{\lVert h\rVert} \\
& \leq\lim_{h\to 0} K\lVert h\rVert\dfrac{\lVert g(A)-g(a)\rVert}{\lVert h \rVert}=0
\end{align*}
This completes the proof.
\end{proof}
\section{Conformable Partial Derivatives}
We begin this section with a correction. In \citep{3}, authors define the conformable derivative of a function with $m$ variables [Definition 2.9]. But the definition is not valid for $\alpha =\frac{1}{2}$ and $x_i<0$. We give the corrected definition of conformable partial derivative of a real valued function with $n$ variables by the following.
\begin{definition}\label{Def4.1}
Let $f$ be a real valued function with $n$ variables and $a=(a_1,...,a_n)\in\mathbb{R}^n$ be a point whose $i$th component is positive. Then the limit
\begin{equation}\label{eq5}
\lim_{\epsilon\to 0}\dfrac{f(a_1,...,a_i+\epsilon a_i^{1-\alpha},...,a_n)-f(a_1,...,a_n)}{\epsilon},
\end{equation}
if it exists, is denoted by $\dfrac{\partial ^{\alpha}}{\partial x^\alpha}f(a)$, and called the $i$th conformable partial derivative of $f$ of order $\alpha\in (0,1]$ at $a$.
\end{definition}
\begin{theorem}\label{Theo4.2}
Let $f$ be a vector valued function with $n$ variables. If $f$ is $\alpha -$differentiable at $a=(a_1,...,a_n)\in\mathbb{R}^n$, each $a_i>0$, then $\dfrac{\partial ^{\alpha}}{\partial x_j^\alpha}f_i(a)$ exists for $1\leq i\leq m$, $1\leq j\leq n$ and the conformable Jacobian of $f$ at $a$ is the $m\times n$ matrix $\Big( \dfrac{\partial ^{\alpha}}{\partial x_j^\alpha}f_i(a)\Big)$.
\end{theorem}
\begin{proof}
Let $f(x_1,...,x_n)=(f_1(x_1,...,x_n),...,f_m(x_1,...,x_n))$. Suppose first that $m=1$, so that $f(x_1,...,x_n)\in\mathbb{R}$. Define $h:\mathbb{R}\longrightarrow\mathbb{R}^n$ by $h(y)=(a_1,...,y,...,a_n)$ with $y$ in the $j$th place. Then $\dfrac{\partial ^{\alpha}}{\partial x_j^\alpha}f_i(a)=D^\alpha (f\circ h)(a_j)$. Hence, by Corollary \ref{Cor3.11},
\begin{align*}
(f\circ h)^\alpha (a_j) & =f^\alpha (h(a_j))\begin{bsmallmatrix}
h_1(a_j)^{\alpha -1} & 0 &  & \hdots &  & 0 \\
\vdots & \vdots & &\ddots & & \vdots \\
0 & 0 & \hdots & h_j(a_j)^{\alpha -1} & \hdots & 0 \\
\vdots & \vdots & &\ddots & & \vdots \\
0 & 0 & & \hdots & & h_n(a_j)^{\alpha -1}
\end{bsmallmatrix}h^\alpha (a_j) \\
& =f^\alpha (a) \begin{bsmallmatrix}
a_1^{\alpha -1} & 0 &  & \hdots &  & 0 \\
\vdots & \vdots & &\ddots & & \vdots \\
0 & 0 & \hdots & a_j^{\alpha -1} & \hdots & 0 \\
\vdots & \vdots & &\ddots & & \vdots \\
0 & 0 & & \hdots & & a_n^{\alpha -1}
\end{bsmallmatrix}\begin{bsmallmatrix}
0 \\
\vdots \\
a_j^{1-\alpha} \\
\vdots \\
0
\end{bsmallmatrix}=f^\alpha (a)\begin{bsmallmatrix}
0 \\
\vdots \\
1 \\
\vdots \\
0
\end{bsmallmatrix}.
\end{align*}
Since $(f\circ h)^\alpha (a_j)$ has the single entry $\dfrac{\partial ^{\alpha}}{\partial x_j^\alpha}f_i(a)$, this shows that $\dfrac{\partial ^{\alpha}}{\partial x_j^\alpha}f_i(a)$ exists and is the $j$th entry of the $1\times n$ matrix $f^\alpha (a)$.

The theorem now follows for arbitrary $m$ since, by Theorem \ref{Theo3.12}, each $f_i$ is $\alpha -$differentiable and the $i$th row of $f^\alpha (a)$ is $(f_i)^\alpha (a)$.
\end{proof}
\section{Conclusion}
In our study, we have extended the idea given by Khalil at al. in \citep{1} and called it multivariable conformable fractional calculus. Multivariable conformable fractional derivative has many interesting properties and is related to classical multivariable calculus. We have given important theorems and corollaries which reveal this relation and we have obtained useful results.


\section*{References}
\begin{thebibliography}{99}
\bibitem{1} R. Khalil, M. Al Horani, A. Yousef, M. Sababheh, A new definition of fractional derivative, J. Comput. Appl. Math 264(2014)65-70.
\bibitem{2} T. Abdeljawad, On conformable fractional calculus, J. Comput. Appl. Math 279(2015)57-66.
\bibitem{3} A. Atangana, D.Baleanu, and A. Alsaedi, New properties of conformable derivative, Open Math. 13(2015)889-898.
\end{thebibliography}
\end{document}